

Let $\Sigma\i P^3$ be a smooth cubic surface.
 It is known that $S$ contains
27 lines. Out of  these lines one can form 36 {\it Schl\"afli double  -
sixes}
i.e., collections $\{l_1,...,l_6\}, \{l'_1,..., l'_6\}$ of 12 lines
such that each $l_i$ meets only $l'_j, \, j\neq i$ and does not
meet $l_j, j\neq i$, see n.0.1 below.
 In 1881 F. Schur proved [S] that  any double - six
gives rise to a certain quadric $Q$ , called {\it Schur quadric}
which is characterized as follows:
for any $i$ the lines $l_i$ and $l'_i$ are
 orthogonal with respect to (the quadratic form defining) $Q$.
\vskip .2cm

\magnification =\magstep1
\baselineskip =13pt

\centerline{\bf SCHUR QUADRICS, CUBIC SURFACES AND RANK 2 VECTOR BUNDLES}
\centerline {\bf OVER  THE PROJECTIVE PLANE}

\vskip 1cm

\centerline {\bf I.Dolgachev and M.Kapranov \footnote{*}{\rm
Both authors were partially supported by NSF.}}

\vskip 1.5cm

Let $\Sigma\i P^3$ be a smooth cubic surface.
 It is known that $S$ contains
27 lines. Out of  these lines one can form 36 {\it Schl\"afli double  -
sixes}
i.e., collections $\{l_1,...,l_6\}, \{l'_1,..., l'_6\}$ of 12 lines
such that each $l_i$ meets only $l'_j, \, j\neq i$ and does not
meet $l_j, j\neq i$, see n.0.1 below.
 In 1881 F. Schur proved [S] that  any double - six
gives rise to a certain quadric $Q$ , called {\it Schur quadric}
which is characterized as follows:
for any $i$ the lines $l_i$ and $l'_i$ are
 orthogonal with respect to (the quadratic form defining) $Q$.
\vskip .2cm

The aim of the present paper is to relate Schur's construction to
the theory of vector bundles on $P^2$ and to generalize this construction
along the lines of the said theory.

\vskip .3cm

 Let us describe the vector bundle interpretation
of the Schur quadric. Note that the first six lines
$\{l_1,...,l_6\}$
of a double - six on $\Sigma$ define a blow-down
$\pi: \Sigma \rightarrow P^2$
which takes the lines $l_i$ into some points $p_i\in P^2$. These
points are in general position i.e. no three of them lie on a line.
Let $\check P^2$ be the dual projective plane and $H_i \i \check P^2$ be
the lines corresponding to $p_i$. The union ${\cal H}$ of these
lines is a divisor with normal crossing in $\check P^2$. Let
$E({\cal H}) = \Omega^1_{\check P_2}(\log {\cal H})$
be the corresponding vector bundle (locally free sheaf)
of logarithmic 1-forms on $\check P^2$.
The twisted bundle
$E = E({\cal H})(-2)$ is a stable
rank 2 bundle on $\check P^2$ with Chern classes $c_1 = -1,
c_2 = 4$ (see [DK]).
For such bundles K.Hulek [Hu 1] has defined the notion of a
{\it jumping line of the second kind} (shortly JLSK).
 This is a line $l\i \check P^2$ such that the restriction of
$E$ to
 the first
infinitesimal neigborhood $l^{(1)}$ of $l$ is not isomorphic to
${\cal O}_{l^{(1)}} \oplus {\cal O}_{l^{(1)}}(-1)$. Hulek has shown that
such lines form a curve $C(E)$ in the projective plane of lines
 in $\check P^2$
i.e. in $P^2$. Now the result is as follows.

\vskip .2cm

\proclaim  Theorem 1. The space $P^3$ containing the cubic surface
$\Sigma$ is
naturally identified with the projectivization of $H^1(\check P^2,
 E(-1))^*$.
Under this identification the Schur quadric $Q$ becomes
dual to the zero locus of the
quadratic form given by the cup-product
$$H^1(\check P^2, E(-1))\otimes H^1(\check P^2, E(-1)) \rightarrow H^2
 \biggl(\check P^2, \bigl(\bigwedge^2E\bigl)(-2)\biggl) =
H^2(\check P^2, {\cal O}(-3)) = {\bf C}.$$
The intersection $\Sigma\cap Q$ is mapped, under the projection
 $\pi:\Sigma\rightarrow P^2$, to the curve of JLSK $C(E)$.

More generally, the whole theory of Hulek [Hu 1] of rank 2 vector bundles
on $P^2$ with odd $c_1$ can be given a "geometric" interpretation involving
some natural generalizations of cubic surfaces, double - sixes and Schur
quadrics. This is done in \S 2 of the paper. This interpretation
implies Theorem 1.

The outline of the paper is as follows. In \S 0 we recall some known
(and less known) facts about cubic surfaces and Schur quadrics. In
\S 1 we give a short overview of Hulek's theory of monads corresponding
to vector bundles with $c_1 = -1$. In \S 2 we give an interpretation
of Hulek's theory mentioned above. In \S 3 we consider bundles of logarithmic
1-forms corresponding to arrangements of $2d$ lines in $P^2$ in general
position. The main result of this section is that all these bundles
satisfy certain condition of
$\Sigma$ - genericity in the sense defined in \S 2, which makes working
with bundles satisfying this condition easier. Finally, in \S 4 we consider
various examples of the previous constructions corresponding to some
special types of vector bundles.

\beginsection

\centerline {\bf \S 0. Cubic surfaces.}

\vskip 1cm

\noindent {\bf 0.1.} Here we recall some standard known facts about
cubic surfaces. All the proofs can be found either in [H], Ch.V, \S 4
or in [M] or can be easily reconstructed by the reader. Let
$p_1,...,p_6$ be six distinct points in the projective plane $P^2$.
 Assume that no
three of these points lie on a line.  Denote by $Z$ the union of the
points $p_i$ and by ${\cal J}_Z\i {\cal O}_{P(V)}$ the
sheaf of ideals of $Z$.
The linear system $P\bigl (H^0({\cal J}_Z(3))\bigl)$
 of cubic curves through $Z$
is of dimension 3 and defines a rational map
$$f: P^2 \rightarrow
P(H^0({\cal J}_Z(3)^*) = P^3$$
 whose image is a cubic surface, denoted $\Sigma$.
The rational map $f$ comes from a regular map $f':
{\rm Bl}_Z(P^2) \rightarrow P^3$
where ${\rm Bl}_Z(P^2)$ is the blow up of $Z$.
Let $\pi: {\rm Bl}_Z(P^2) \rightarrow P^2$ be the projection.
If we further
assume that the points $p_i$ do not lie on a conic then $f'$
 is an isomorphism
and $\Sigma$ is nonsingular. If $p_i$ do lie on a conic then $\Sigma$
is singular
and $f'$ blows down this conic to a singular point of $\Sigma$.

\vskip .3cm

Suppose $\Sigma$ is nonsingular. Then $\Sigma$ has 27 lines on it.
 They can be grouped
into three subsets:
$$\{l_1,...,l_6\},\,\, \{l'_1,...,l'_6\},\,\, \{m_{ij}, \,
 1\leq i<j\leq 6\}.\eqno (0.1)$$
The lines $l_i$ are the images under $f'$ of the exceptional lines $\pi^{-1}
(p_i)$. The lines $l'_i$ are images under $f'$ of proper transforms of
the conics $C_i\i P^2$ passing through $Z-\{p_i\}$.  Finally the lines
 $m_{ij}$ are images of the proper transforms of the lines $<p_i, p_j>$
joining the points $p_i$ and $p_j$.

The first two groups of lines form a {\it double - six} which means that
$$l_j\cap l_j = \emptyset, \quad
\l'_i \cap \l'_j = \emptyset,\quad l_i\cap l'_j \neq \emptyset
\quad\quad {\rm iff}
\quad i\neq j.\eqno (0.2)$$
Every set of 6 disjoint lines on $\Sigma$ can be included in
 a unique double -
six from which
$\Sigma$ can be reconstructed uniquely.
 There are 36 double - sixes of $\Sigma$. Every double - six defines
two regular birational maps $\pi_1: \Sigma \rightarrow P^2$, $\pi_2: \Sigma
\rightarrow P^2$, each blowing down one of the two sixes (sixtuples
of disjoint lines) of the double - six. The birational map
$\pi_2 \circ \pi_1^{-1}: P^2 \rightarrow P^2$ is given by the linear
system of quintics with double points at $p_i$.
The two collections of 6 points in $P^2$ given by $\{\pi_1(l_i)\}$
and  $\{\pi_2 (l'_i)\}$ are associated to each other
 in the sense of Coble (cf.[DO] [DK]).

\vskip .3cm

\noindent {\bf 0.2.} Here we shall discuss somewhat less known facts about
the determinantal representation of a cubic surface [B]. A modern treatment
of this can be found in [G] [Gi]. Consider the
homogeneous ideal of the subscheme $Z$ i.e.
$$I_Z = \bigoplus_{n\geq 0} H^0(P^2, {\cal J}_Z (n))\eqno (0.3)$$
in the graded ring $R = {\bf C}\,\,[T_0, T_1, T_2]$.
 It is easy to see that
the ring $R/I_Z$ is Cohen - Macaulay hence of homological
dimension 1.
Any four linearly independent cubic forms vanishing on $Z$
represent a minimal
set of generators of $I_Z$. According to the Hilbert-Burch
 theorem, see [Kap],
 the ideal $I_Z$ is generated by the
 maximal minors
of some $3\times 4$ matrix of homogeneous linear forms.
 In other words, we
have a resolution
$$0\rightarrow R(-4)^3 \rightarrow R(-3)^4 \rightarrow I_Z
 \rightarrow 0.$$
This resolution gives the resolution of the  sheaf
${\cal J}_Z(3)$:
$$0\rightarrow {\cal O}_{P(V)}(-1)^3 \rightarrow {\cal O}_{P(V)}^4
\rightarrow {\cal J}_Z(3) \rightarrow 0.$$
We can rewrite this resolution in the form
$$0\rightarrow {\cal O}_{P^2}(-1)\otimes I^* \buildrel
\gamma\over\rightarrow  {\cal O}_{P^2} \otimes L^*
 \rightarrow {\cal J}_Z(3)
\rightarrow 0 \eqno (0.4)$$
where vector spaces $I^*$ and $L^*$ of respective
dimensions 3 and 4 are
defined intrinsically as follows:
$$L^* = H^0(P^2, {\cal J}_Z(3)); \eqno (0.5)$$
$$I^* = {\rm Ker} \{H^0(P^2, {\cal O}(1)\otimes L^*) \rightarrow H^0(P^2,
{\cal J}_Z (4)) \}. \eqno (0.6)$$

It will be convenient for us to regard henceforth our projective
plane $P^2$ as $P(V^*)$ where $V$ is a 3-dimensional vector space.
With this choice of notation, the map $\gamma$ in (0.4) is given by a
linear map $I^*\otimes V^* \rightarrow L^*$. We shall be more
interested in the transpose of this map which we denote by
$$g: L \longrightarrow I \otimes V = {\rm Hom} (V^*, I). \eqno (0.7)$$
Choosing bases in $V, I$ we can regard $g$ as a 3 by 3 matrix of linear
forms on $L$. Here is the classical result on the determinantal
representation.
\vskip .3cm

\proclaim 0.3.~Proposition. The map $g$ is an embedding. The locus
$$\Sigma = \{x\in P(L): \,\,\, {\rm rank}\,\,\, g(x) \leq 2\}
\eqno (0.8)$$
is a nonsingular cubic surface in $P(L) = P^3$ isomorphic to
 ${\rm Bl}_Z(P(V^*)$.
An explicit blow-down $\pi_1: \Sigma \rightarrow P(V^*)$
 takes $x\in \Sigma$
into {\rm Ker} $g(x) \in P(V^*)$. It is isomorphism
outside the set $Z = \{p_1,
.., p_6\}\i P(V^*) = P^2$ (see n. 0.1). The dual blow-down
  $\pi_2:
\Sigma \rightarrow P(I^*)$ takes $x\in \Sigma$ into
$({\rm Im} g(x))^\bot \in P(I^*)$. It is an isomorphism
outside a six - element set $Z^{as} = \{q_1, ..., q_6\}
\i P(I^*)$ (this is the set associated to $Z$).

Note that a given cubic surface $\Sigma\i P^3$ has many
non-equivalent
determinantal representations corresponding to different ways of blowing
down $\Sigma$ onto a $P^2$ (i.e. to different choices of a double - six).

\vskip .3cm

\noindent {\bf 0.4.} All the other attributes of the cubic
surface $\Sigma$
can be easily found from the map $g$. For example,
the set $Z$ can be recovered
in terms of $g$ as follows. Consider the partial transpose
of (0.7):
$$g_V: V^* \rightarrow I\otimes L^* = {\rm Hom}(L, I).$$
Then
$$Z = \{ z\in P(V^*): \,\, {\rm rank}\,\,\, g_V(z)\leq 2\}.
\eqno (0.9)$$
The 12 lines of the double - six can be written in the form
$ A_z = P({\bf A}_z), A'_z = P({\bf A}'_z), \, z\in Z$
where ${\bf A}_z$ and ${\bf A}'_z$ are 2-dimensional
 vector
subspaces in $L$ defined for $z\in Z$ as follows:
$${\bf A}_z = \{x\in L: \, g(x) \in ({\rm Im}\,\,\, g_V(z))\otimes V\};
\eqno (0.10)$$
$${\bf A}'_z = \{x\in L: \, g(x) \in I\otimes z^\bot\},\eqno (0.11)$$
where $z^\bot\i V$ is the 2-plane orthogonal to $z\in P(V^*)$. Thus if
$Z = \{p_1,...,p_6\}$ then the line $A_{p_i}$ is what was
 denoted in n.0.1
by $l_i$ and the line $A'_{p_i}$ is $l'_i$.

\vskip .2cm

The classical theorem of F. Schur [S] can be stated as follows.

\proclaim 0.5.~Theorem. There exists a unique, up to a
scalar factor, symmetric
bilinear form $C(x,y)$ on $L$ with the following property:
$C(x,y)=0$ whenever
$x\in A_z, y\in A'_z$ for some $z\in Z$ (i.e. the corresponding
lines
of the double - six are orthogonal with respect to $C$).
This form
is non-degenerate.

\noindent {\sl Proof.} a) Non-degeneracy: Suppose
such a form $C$ exists and is
degenerate. Let $K$ be the kernel of $C$. Suppose
 ${\rm dim}\,\,\, K = 1$. Then for any
2-dimensional subspace $\Lambda\i L$ not meeting $K$
its orthogonal (with respect to $C$) is a 2-subspace
containing $K$. Since $P(A_z), \, P(A'_z)$ form a
double - six, $K$ can lie on no more than one among
the $A_z$ and no more than one among the $A'_z$. Hence
 there is a 4-element subset $Z_0 \i Z$ such that for
$z\in Z_0$ both $A_z$ and $A'_z$ do not contain $K$.
For such $z$ the space
$A'_z$ should coincide with $A_z^\bot$ and hence contain
$K$. Hence for
$z_1\neq z_2\in Z_0$ we have $A'_{z_1} \cap A'_{z_2}
\neq \{0\}$ which is
a contradiction. The cases dim $K$ = 2,3 are similar and left to the reader.

\vskip .1cm

b) Uniqueness: If there are two non-proportional forms $C_1, C_2$ with the
required property then for any $\lambda, \mu$ the linear combination
$\lambda C_1 + \mu C_2$ also satisfies this property. However, there will
be always such $\lambda, \mu$ that the linear combination is non-zero
but degenerate. This contradicts a).

\vskip .3cm

\noindent {\bf 0.6.} It remains to prove the existence part of Theorem 0.5.
To do this, let us take the second symmetric power of
the map $g$ in (0.7) and use the natural decomposition
$$S^2(I\otimes V) = \biggl(\bigwedge^2 I \otimes \bigwedge^2 V\biggl)
\quad \oplus \quad
\bigl(S^2 I \otimes S^2 V\bigl).\eqno (0.12)$$
By projecting $S^2g$ to the first summand, we get a linear map
$$S^2L \longrightarrow S^2(I\otimes V) \longrightarrow
\biggl(\bigwedge^2 I \otimes \bigwedge^2 V\biggl).\eqno (0.13)$$
Note that dim $S^2L = 10$, and dim$(\bigwedge^2 I
 \otimes \bigwedge^2 V) = 9$.
Hence the map (0.13) has  non-trivial kernel. (We shall see later that
this kernel is in fact 1-dimensional).

\proclaim 0.7.~Proposition. If $B$ is a non-zero form from the kernel of
(0.13) then $B: L^* \rightarrow L$ is invertible and $C=B^{-1}:
L\rightarrow
L^*$ is a bilinear form on $L$ satisfying the conditions of Theorem 0.5.

We shall concentrate on the proof of this proposition.

\vskip .3cm

\noindent {\bf 0.8.} A form $B\in S^2L$ lying in the kernel of (0.13)
is classically called "apolar to all the quadratic forms given
 by $2\times 2$
minors of $g$", cf. [B]. In general, if $E$ is a vector space then
quadratic forms $G\in S^2E,\,\, H\in S^2 E^*$ are called {\it
apolar} if $(G,H)_2 = 0$ where $(\cdot , \cdot)_2$ is the natural pairing
$S^2E \otimes S^2E^* \rightarrow {\bf C}$.
Note the particular case when $G$ has rank 2 i.e. $G = e\cdot f$ is the
symmetric
product of two vectors $e,f\in E$. In this case the apolarity of
 $G$ and $H$
means that $H(e,f) = 0$.

We shall need a different description of the map dual to (0.13).
Let us denote
this map by
$$\delta: \bigwedge^2 I^* \otimes \bigwedge^2 V^*
\longrightarrow S^2L^*.
\eqno (0.14)$$
Let us chose volume forms on $V$ and $I$. Then we can write
$\bigwedge^2 V^* =
V,\,\, \bigwedge^2 I^* = I$. It is immediate
to see that there are identifications
$$\bigwedge^2 V^* \cong V \quad =
\quad H^0\bigl(P(V^*), {\cal O}_{P(V^*
)}
 (1)\bigl);\eqno (0.15)$$
$$\bigwedge^2_{} I^* \cong I \quad
= \quad H^0\bigl(P(V^*), {\cal J}_Z^2(5)\bigl);
\eqno (0.16)$$
$$S^2L^* \quad = \quad H^0\bigl(P(L),
 {\cal O}(2)\bigl) \quad \cong \quad H^0\bigl(P(V^*),
{\cal J}_Z^2 (6)\bigl).\eqno
(0.17)$$
 Indeed, (0.15) follows by definition of ${\cal O}(1)$;
the identification (0.16) expresses the fact that the Cremona
transformation $\pi_2 \circ \pi_1^{-1}: P(V^*) \rightarrow P(I^*)$
is given by the linear system of quintics with singular points $p_i$,
see n. 0.1. Finally,  to see (0.17) we note that the embedding
of the cubic surface $\Sigma$ into $P(L) = P^3$ is given by
the linear system of cubics in $P(V^*)$ through $p_i$,
so $L^*$ is the space of cubic polynomials on $V^*$
vanishing at $p_i$. The second symmetric power of this space
maps therefore to the space of polynomials of degree 6
vanishing at $p_i$ together with their first derivatives
i.e, to the RHS of (0.17);
this map is easily seen to be an isomorphism.

\proclaim 0.9.~Lemma. Under identifications (0.15) - (0.17) the map
$\delta$ corresponds to the multiplication map
$$H^0 \bigl(P(V, {\cal O}(1)\bigl) \otimes H^0
\bigl(P(V^*), {\cal J}_Z^2(5)\bigl)
\rightarrow H^0 \bigl(P(V^*), {\cal J}_Z^2 (6)\bigl).$$

In other words, quadrics in $P(L)$ are identified
 with sextics in $P(V^*)$
with double points at $p_i\in Z$ and quadrics from
the image of $\delta$
correspond to sextics containing a line.

\vskip .3cm

\noindent {\sl Proof of the lemma:} We have the
commutative diagram
$$\matrix{&L&\buildrel {\rm Sq} \over\longrightarrow & S^2L \cr
g& \big\downarrow &&\big\downarrow\cr
&{\rm Hom} (V^*, I) &\buildrel \lambda \over\longrightarrow
& \bigwedge^2 {\rm Hom} (V^*, I)}\eqno (0.18)$$
where the map Sq takes $x\mapsto x^2$, the map $\lambda$
takes $\phi \mapsto \bigwedge^2 \phi$
 and the map on the right is the same
as in (0.13). We keep the volume forms in $I$ and $V$
and identify correspondingly the spaces $\bigwedge^2I$ and $\bigwedge
^2 V$ with $V^*$ and $I^*$. For any $\phi \in {\rm Hom}(V^*, I)_2$
the second exterior power $\bigwedge^2 \phi \in I^* \otimes V^*$
is a tensor of rank 1. Hence it can be written in the form
$i^* \otimes v^*$ for some $i^* \in I^*, \,\, v^* \in V^*$.
This shows that the restriction of the map $\lambda \circ g$
to the cubic surface $\Sigma \i P(L)$ coincides with the
composition
$$\Sigma \buildrel \pi_1\times \pi_2 \over \longrightarrow
P(I^*) \times P(V^*) \buildrel {\rm Segre} \over
\longrightarrow P(I^* \otimes V^*)\eqno (0.19)$$
where $\pi_j$ are the blow-downs from n. 0.3.

The map $\lambda \circ g: P(L) \rightarrow P(I^*\otimes V^*)$
 is given by the linear system ${\cal Q}$
of quadrics which is the
projectivization of the image of the linear  map $\delta$
from (0.14).
The system ${\cal Q}$ is spanned by the $2\times 2$
 minors of the matrix of linear forms on $L$ defining the
determinantal representation of $\Sigma$.
In other words, the preimage of the linear system of
hyperplane sections
of $P(I^*\otimes V^*)$ under $\lambda\circ g$ is
the linear system of quadric sections on $\Sigma$
which is (the projectivization of)
the image of the canonical pairing
$$H^0 \bigl(P(I^*), {\cal O}(1)\bigl)
\otimes H^0 \bigl(P(V^*), {\cal O}(1) \bigl)\quad
\longrightarrow  \quad H^0\bigl(\Sigma, {\cal O}(2)\bigl).$$
By Theorem 0.3, we can make an identification of the projective
spaces $P(I^*)$ and
\hfill\break $P\bigl(H^0\bigl(P(V^*)$, ${\cal J}_Z^2(5)
\bigl)^*\bigl)$.
Under the rational map $P(V^*) \rightarrow P(I^*)$ given
by the linear system of sections of ${\cal J}_Z^2(5)$,
 zeroes of these sections are preimages of the lines in $P(I^*)$
and the resulting map
$$H^0\bigl(P(V^*), {\cal O}(1)\bigl) \otimes H^0
\bigl(P(V^*), {\cal J}_Z^2(5)\bigl)
\rightarrow H^0 \bigl(\Sigma, {\cal O}(2)\bigl) =
H^0\bigl(P(V^*), {\cal J}_Z^2(6)\bigl)$$
coincides with the natural multiplication map from the assertion of the
lemma. So the lemma is proven.

\vskip .3cm

\noindent {\bf 0.10.} We continue to
prove Theorem 0.5 and shall now use  Lemma 0.9.
Let us consider some particular sextics with
double points at $Z = \{p_1,...,p_6\}$. Let $C_i$ be the unique conic
through $Z-\{p_i\}$. We can take a sextic curve which is the union of
two lines $<p_i, p_j>,\quad <p_k,p_s>$ and  two conics $C_k, C_s$. By
means of (0.17) this sextic corresponds to some quadric $Q_{ij, ks}$.
Moreover, since the quintic $<p_k, p_s> \,\cup \,\,C_k \cup
\,\, C_s$
belongs to the linear system of quintics singular at points of $Z$,
 the quadric
$Q_{ij, ks}$ lies in the image of  the map $\delta$
from (0.14). Now let us take $j=s$.
Then our sextic can be represented as the union of two cubic curves
through $Z$ namely
$$<p_i, p_j> \cup \,\,\,C_j \quad {\rm and} \quad <p_k, p_j> \cup
\,\,\,C_k.$$
Since such cubics correspond to hyperplanes in $P(L)$, we conclude
that the quadric $Q_{ij, kj}$ is in fact the union of two planes,
say $H_{ij}$ and $H_{kj}$. Moreover, $H_{ij}$ cuts out the cubic surface
$\Sigma$ along 3 lines $l_i, l'_j, m_{ij}$ (see n. 0.1). The plane
$H_{kj}$ cuts out the lines $l_j, l'_k, m_{kj}$ on $\Sigma$.
Since the (quadratic form defining the) quadric $Q_{ij, kj} = H_{ij} \cup
H_{kj}$ is apolar to our chosen $B\in S^2L$, we conclude that the
equations of $H_{ij}$ and $H_{kj}$ (belonging to $L^*$) are $B$ -
orthogonal.

\vskip .3cm

\noindent {\bf 0.11.} Let us now prove Proposition 0.7 and hence
Theorem 0.5. The form $B$ is a linear map $L^* \rightarrow L$. For
any linear subspace $U\i L$ we define its {\it polar subspace}
(with respect to $B$) to be
$$ U^\bot = \sum_{\xi\in L^*, \xi|_U = 0} {\bf C}\cdot B(\xi).$$
If $B$ is non-degenerate then $U^\bot$ is the usual orthogonal
complement of $U$ with respect to the inverse form $B^{-1} \in
S^2L^*$. If $B$ is degenerate and $K\i L^*$ is its kernel then every
$U^\bot$ is contained in $K^\bot$ where $K^\bot\i L$ is the orthogonal
to $K$ with respect to the usual pairing of $L^*$ and $L$.

We shall use the same terminology and notation for projective
subspaces in $P(L)$. In particular, if $H\i P(L)$ is a hyperplane
whose equation does not lie in $K = {\rm Ker}\,\, B$ then $H^\bot$
is a point called the {\it pole} of $H$.

Let us prove that $B$ is non-degenerate. Let the double - six be
$\{l_1,...,l_6\}, \,\, \{l'_1,...,l'_6\}$. Consider the plane
$H_{12}$ spanned by lines $l_1$ and $l'_2$ (which intersect).
Its equation (in $V^*$) is orthogonal with respect to $B$ to equations
of similar planes $H_{21}, H_{23}, H_{31}$ (see n. 0.10). Hence
$H_{12}^\bot = H_{21}\cap H_{23}\cap H_{31}$ and this intersection
is easily seen to be the point $l'_1\cap l_2$. In this way we show that
each $l'_i\cap l_j$ is the pole of some plane. Since the points
$l'_i\cap l_j$ span $P(L)$, the form $B$ cannot be degenerate (see
above).

It remains to show that $l_i^\bot = l'_i$. We have already seen that the
point $l'_1\cap l_2$ is the pole of the plane $H_{12}$ spanned by
$l_1\cap l'_2$. Similarly, $l'_1\cap l_3$ is the pole of $H_{13} =
{\rm Span} ( l_1, l'_3)$ Hence $l'_1 = {\rm Span}
(l'_1\cap l_2, l'_1\cap l_3)$ is the orthogonal complement of
$H_{12}\cap H_{23} = l_1$. Similarly we prove that $l'_i = l_i^\bot$
for other $i$.

\vskip .2cm

Theorem 0.5 is completely proven.
The reader should compare this rather cumbersome proof with a more
straightforward one based on the theory of vector bundles
(Theorem 2.17 below).

\proclaim 0.12.~Definition. The quadric $Q\i P(L)$ defined by $C(x,x)=0$
where $C$ is the quadratic form given by Theorem 0.5, is called the Schur
quadric (associated with the double - six $\{A_z, A'_z\}$).

\noindent {\bf 0.13. Example.}
Let us consider the following 4-dimensional
space $L$:
$$L = \biggl\{(x_1,...,x_5) \in {\bf C}^5:  \quad
\sum x_i =0 \biggl\}$$
and define the cubic surface $\Sigma\i P(L)$ by the equation $\sum x_i^3=0$
(the Clebsch diagonal surface).
The symmetric group $S_5$ acts on ${\bf C}^5$ by permutations of
coordinates and preserves $L$ and $\Sigma$.
The  line
$$l = \biggl\{x\in P(L): \quad x_1 + {1 + \sqrt 5 \over 2} x_2 + x_3 =
 x_2 + {1 + \sqrt 5 \over 2}x_3 + x_4 = 0\biggl\} $$
 lies on $\Sigma$ and so do all the lines obtained from $l$
by the action of $S_5$. It is known [Bu] that the $S_5$ -
orbit of $l$ consists of 12 lines which  form a
double - six. Their equations can be found in [B], p.168. The
two sextuples of lines constituting this double - six are
orbits of the alternating group $A_5 \i S_5$. So one sextuple
is the $A_5$ - orbit of $l$ and the other is the $A_5$ - orbit
of the line
$$l' = \biggl\{ x: \,\,\,\, x_1 + x_2 + {1 - \sqrt 5 \over 2}x_4
= {1 - \sqrt 5 \over 2}x_1 + x_3 + x_4 = 0\biggl\}$$
So $l'$ is  line of the second sextuple corresponding to
$l$ (because $l\cap l' = \emptyset$). The lines $l$
and $l'$ are orthogonal with respect to the bilinear form
$C(x,y) = \sum_{i=1}^5 x_iy_i$
on $L$. By symmetry, all the other corresponding pairs
of lines of our double  - six are  also orthogonal with respect to $C$.
Thus the Schur quadric $Q$ associated to this
double - six is given by the equation
$\sum_{i=1}^5 x_i^2 = 0$.

\beginsection

\centerline {\bf \S 1. An overview of Hulek's theory.}

\vskip 1cm

\noindent {\bf 1.1.} Let $E$ be a
stable rank 2 vector bundle on $P^2 = P(V)$ with
$c_1(E) = -1, \,\, c_2(E) = n$.
According to Le Potier [L] and Hulek [Hu 1],
the bundle $E$ can be realized as the middle cohomology of a monad
$$H\otimes {\cal O}_{P(V)}(-1) \buildrel \alpha
\over \longrightarrow
M\otimes \Omega^1_{P(V)}(1) \buildrel \beta \over
\longrightarrow H'\otimes {\cal O}_{P(V)}.\eqno (1.1)$$
where
$$H = H^1 \bigl(E(-2)\bigl) \cong {\bf C}^{n-1} ,\quad
 M = H^1 \bigl(E(-1)\bigl) \cong {\bf C}^n, \quad H' =
H^1(E)\cong {\bf C}^{n-1}\eqno (1.2)$$
and the maps $\alpha$ and $\beta$ are defined as follows.
Let $\Omega^1(1)$ be identified with $\Theta (-2)$
where $\Theta$ is the tangent bundle of $P(V)$.
Let $t: V \otimes {\cal O}_{P(V)}(-1) \rightarrow
\Omega^1(1)$ be the Euler homomorphism twisted by
${\cal O}(-1)$ (see [OSS]).
It allows one to identify
$${\rm Hom} \bigl(H\otimes {\cal O}_{P(V)}(-1), M \otimes \Omega^1(1)
\bigl)
\quad \cong\quad {\rm Hom}_{\bf C} (H\otimes V^*, M).\eqno (1.3)$$
The map $\alpha$ is induced by the cup - product
$$a: H^1\bigl(E(-2)\bigl) \otimes V^* \,\, = \,\,
H^1\bigl(E(-2)\bigl) \otimes H^0
\bigl({\cal O}_{P(V)}(1)\bigl) \longrightarrow
H^1\bigl(E(-1)\bigl). \eqno (1.4)$$

\noindent
Similarly, we have a map $t^*: \Omega^1(1) \rightarrow
V^* \otimes {\cal O}_{P(V)}$ which allows us to identify
$${\rm Hom} \bigl(M \otimes \Omega^1(1), H' \otimes
{\cal O}_{P(V)}\bigl)
\quad \cong \quad {\rm Hom}_{\bf C}(M\otimes V^*, H').
\eqno (1.5)$$
After this identification the map $\beta$ is induced by the
cup - product
$$b: H^1\bigl(E(-1)\bigl) \otimes V^* \,\, = \,\,
H^1\bigl(E(-1)\bigl) \otimes H^0
\bigl({\cal O}_{P(V)}(1)\bigl) \longrightarrow
H^1(E). \eqno (1.6)$$
The cup - product pairing
$$B: M \otimes M \,\, = \,\, H^1\bigl(E(-1)\bigl) \otimes
H^1\bigl(E(-1)\bigl)
\longrightarrow H^2 \biggl(\biggl(\bigwedge^2 \,\,E\biggl)
(-2)\biggl) \,\,\, = $$
$$ = \,\,\, H^2\bigl({\cal O}_{P(V)}(-3)\bigl) \quad  =
\quad {\bf C}
\eqno (1.7)$$
is a symmetric non - degenerate bilinear form on $M$. We
regard it as an isomorphism
$$B: M \rightarrow M^*. \eqno (1.8)$$
The spaces $H$ and $H'$ are dual to each other by means of the
Serre duality and the isomorphism
$$E \quad \cong \quad E^* \otimes \bigwedge^2 E \quad \cong
\quad E^*(-1).$$
With respect to the constructed pairings the monad (1.1) is self -
dual in the sense that $\beta = \alpha^* (-1)$. Equivalently,
if $\lambda\in V^*$ and $a(\lambda): H \rightarrow M$ is the
linear map defined by the pairing $a$ and similarly
 $b(\lambda): M \rightarrow H'$  is the map defined by $b$
then
$$b(\lambda) = a(\lambda)^* \circ B.$$
This shows that the monad (1.1) is completely determined by the
pairing (1.4) and the symmetric bilinear form $B$. The pairing
must satisfy the following properties (cf. [Hu 1]):

\vskip .2cm

\item{$(\alpha 1)$} The map $a(\lambda)$ is injective
for generic $\lambda \in V^*$.
\item{$(\alpha 2)$} For any $h \in H$ the map $a_H(h):
V^* \rightarrow M$ defined by the pairing $a$ is of rank $\geq 2$.
\item{$(\alpha 3)$} For any $\lambda, \lambda' \in V^*$ we have
$b(\lambda') \circ a(\lambda) = b(\lambda) \circ a(\lambda')$
where $b(\lambda) = a(\lambda)^* \circ B$ and similarly for
$b(\lambda')$.

\proclaim 1.2.~Theorem. Let $V, H, M$ be linear spaces of
respective dimensions $3, n-1$ and $n$ and $n\geq 2$. Let
us fix a non-degenerate symmetric bilinear form $B$ on $M$.
By assigning to each $a \in {\rm Hom}(H\otimes V^*, M)$
satisfying $(\alpha 1)$ -- $(\alpha 3)$ the map
$$\alpha = \bigl({\rm Id} \otimes t\bigl) \circ
 \bigl(a\otimes {\rm Id}
\bigl): \,\, \, H \otimes {\cal O}_{P(V)}
(-1) \longrightarrow M \otimes \Omega^1 (1),$$
we get a bijective correspondence between equivalence classes
of self - dual monads (1.1) modulo action of the group
$O(M,B) \times GL(H)$ and isomorphism classes of stable rank 2
vectpr bundles $E$ on $P(V)$ with $c_1(E) = -1$ and
$c_2(E) = n$.

\vskip .2cm

\noindent {\bf 1.3.} Let $l$ be a line in $P^2$ and $E$ be a stable
bundle as in Theorem 1.2.
Let $\lambda \in V^*$ be a linear form defining $l$.
 We have a canonical exact
sequence
$$0 \longrightarrow E(-1)
\buildrel \lambda\over\longrightarrow E \longrightarrow
E|_l \longrightarrow 0,$$
which together with the fact $H^0(E) = 0$ which follows from
the stability of $E$, gives an isomorphism
$$H^0\bigl( E|_l \bigl) =
{\rm Ker}\bigl\{ H^1(E(-1)) \rightarrow H^1(E)\bigl\} =
{\rm Ker} \bigl\{ a(\lambda): M \rightarrow H\bigl\}.
\eqno (1.9)$$
Since $E|_l \cong {\cal O}(p) \oplus {\cal O}(q)$ with $p+q = -1$,
we obtain that
$$E|_l \cong {\cal O} \oplus {\cal O}(-1) \quad \quad
\Longleftrightarrow \quad\quad {\rm rank} a(\lambda) = n-1.$$
A line $l$ is called a {\it jumping line} if $E|_l \neq
{\cal O} \oplus {\cal O}(-1)$. It follows from the Grauert -
M\"ulich theorem [OSS] that the set of jumping lines is a proper
Zariski closed subset of the dual plane $P(V^*)$. This set is
known to be 0 -dimensional for a generic $E$.

\vskip .3cm

\noindent {\bf 1.4.} In [Hu 1] the notion of a
{\it jumping line of the second kind} (shortly JLSK)
was introduced. Let $l^{(1)}$ be the first infinitesimal
neighborhood of $l$ in $P(V)$. We use the exact sequence
$$0\rightarrow {\cal O}_{P(V)}(-2) \buildrel \lambda^2
\over\longrightarrow {\cal O}_{P(V)} \rightarrow
{\cal O}_{l^{(1)}} \rightarrow 0\eqno (1.10)$$
to obtain
$$H^0\bigl(E|_{l^{(2)}}\bigl)  = {\rm Ker}\bigl\{s(\lambda): H^1(E(-2))
\rightarrow H^1(E)\bigl\}.\eqno (1.11)$$
Here the map $s(\lambda)$ corresponds to the canonical
pairing
$$S^2(V^*) \otimes H^1 \bigl(E(-2)\bigl)\quad
 \longrightarrow  \quad H^1(E))$$
evaluated at $\lambda^2$. In the notation of the
previous subsections, $s(\lambda)$ is the composition
$$a(\lambda)^* \circ B \circ a(\lambda): \,\,\,
H\rightarrow M \rightarrow M^* \rightarrow H^*.$$

We say that $l$ is a JLSK if $s(\lambda)$ is not bijective.
Since the source and target of $s(\lambda)$ have the
same dimension, $l$ is a JLSK if and only if
$H^0 \bigl(E|_{l^{(1)}}\bigl) \neq 0$.

\vskip .2cm

Let us introduce a rational quadratic map
$$\gamma: P(V^*) \rightarrow P(S^2H^*),\quad \lambda \mapsto s(\lambda)$$
By property $(\alpha 1)$, for a generic line $l\in P(V^*)$
the value $\gamma (\lambda)$ is well defined and
is an non-degenerate quadric in $P(H)$.
We denote by $C(E)$ the set of all JLSK of $E$.
Thus outside a finite set of points in $P(V^*)$ the
set $C(E)$ is equal to the preimage, under $\gamma$,
of the locus of degenerate quadrics in $P(H)$. So
we get that $C(E)$ is a closed subscheme in $P(V^*)$ defined
by the equation $\det \gamma (l) = 0$. We shall consider $C(E)$
as a closed subscheme of $P(V^*)$ defined by this equation.
So $C(E)$ is a (possibly reducible) curve of degree $2n-2$
containing the set of jumping lines of $E$ in the usual sense.

\vskip .3cm

\noindent {\bf 1.5.} One can give another interpretation of the curve
$C(E)$. Consider the rational map
$$\sigma: P(V^*) \rightarrow P(M^*),\quad \lambda \mapsto {\rm Im}
(a(\lambda))^\bot \i M^*.$$
It is defined on the complement of the set of jumping lines of $E$.
A non-jumping line $l$ is a JLSK of and only if the hyperplane
$\sigma(l) \i P(M)$ is tangent to the quadric
defined by
$B(m,m) = 0$. Let us denote by $Q$ the dual quadric in $P(M^*)$
(which parametrizes the hyperplane tangent to $\{B(m,m) = 0\}$;
so it is given by the inverse quadratic form $C= B^{-1}$).
Then

\vskip .1cm

\centerline { $l$ is a JLSK $\quad$ if and only if
$\quad$  $\sigma (l) \in Q$.}

\beginsection

\centerline {\bf \S 2. Generalized Schur quadrics and cubic surfaces.}

\vskip 1cm

\noindent {\bf 2.1.} Let $E$ be a stable rank 2 vector bundle
 on $P^2 = P(V)$ with $c_1 = -1,\,
c_2 = n$. As we mentioned in the previous section, its
 monad (1.1) defines (and is uniquely defined by) the following
linear algebra data: a linear map (tensor)
$$a: H\otimes V^* \rightarrow M \eqno (2.1)$$
and a quadratic form (the cup - product)
$$B: M\otimes M \rightarrow {\bf C}. \eqno (2.2)$$
Our aim in this section is the study of the geometry of some algebraic
varieties naturally associated to $a$ and $B$ (and hence to $E$).

\vskip .3cm

\noindent {\bf 2.2.} We denote by $Q\i P(M^*)$ the quadric
 defined by the equation $C(m,m) = 0$ where $C$ is the quadratic
 form on $M^*$ inverse to $B$, see n.1.5. We
shall call $Q$
 the {\it Schur quadric} of $E$. We shall see later in this
section
how the classical Schur quadric of a double - six is a particular case of
this construction.

\vskip .3cm

\noindent {\bf 2.3.} By taking various partial transposes of the tensor
$a$, we construct the following linear operators:
$$a_M: M^* \rightarrow H^* \otimes V = {\rm Hom} (H,V); \eqno (2.3)$$
$$a_V: V^* \rightarrow H^*\otimes M = {\rm Hom} (H,M); \eqno (2.4)$$
$$a_H: H \rightarrow M\otimes V = {\rm Hom}(M^*, V). \eqno (2.5)$$
These operators define determinantal varieties in $P(M^*), P(V^*),
P(H)$ consisting of points whose images (under the corresponding $a$)
are operators not of maximal rank. Before going into details, let us recall
some well known facts about varieties of matrices of given rank.

Let $L_1, L_2$ be vector spaces of respective dimensions $n_1, n_2$. We
denote by \hfill\break ${\rm Hom} (L_1, L_2)_r
\i {\rm Hom}(L_1, L_2)$ the variety of
linear maps of rank $\leq r$. We assume that
 $r \leq {\rm min} (n_1, n_2)$.
Then the following is true [ACGH] [R].

\proclaim 2.4.~Proposition. \item{a)} The codimension of
${\rm Hom} (L_1, L_2)_r$  in ${\rm Hom}(L_1, L_2)$ is equal to $(n_1-r)
(n_2-r)$.
\item{b)} ${\rm Hom} (L_1, L_2)_r$ is irreducible and Cohen - Macaulay;
\item{c)} The degree of (the projectivization of)
${\rm Hom} (L_1, L_2)_r$ is equal to
$$\prod_{i=0}^{n_1-r-1} { (n_2+i)! \,\,i!\over (r+i)!\,\, (n_2-r-i)!}.$$
\item{d)} Let $\phi \in {\rm Hom} (L_1, L_2)_r$ be a
linear map of
rank $k\leq r$. Then the multiplicity of ${\rm Hom}(L_1, L_2)_r$
 at $\phi$
is given by
$${\rm mult}_\phi ({\rm Hom}(L_1, L_2)_r) = \prod_{i=0}^{n_1-r-1}
{(n_2-k-i)! \,\,\, i! \over (r-k-1)! (n_2-r-i)!}.$$

\noindent {\bf 2.5.} Let us return to the situation of n. 2.3. We
define the variety $\Sigma \i P(M^*)$ as follows
$$ \Sigma = \{\mu\in P(M^*): {\rm rank} \,\, a_M(\mu) \leq 2\}.
 \eqno (2.6)$$
This is an analog of a cubic surface in $P^3$, cf. Proposition 0.3.

Note that ${\rm dim}\,\, M = n,\,\, {\rm dim} \,\,H = n-1,\,\, {\rm dim}
\,\, V = 3$. Therefore, by Proposition 2.4, the variety ${\rm Hom}(H,V)_2$
has codimension $n-3$ in ${\rm Hom} (H,V)$ and so dim $\Sigma \geq 2$.
Generically, one would expect that dim $\Sigma = 2$.

We shall call the tensor $a$ (and the bundle $E$) $\Sigma$ {\it - generic}
if for any $\mu \in \Sigma$ the rank of $a_M(\mu)$ is exactly 2.
We shall see in section 3 that if $n$ is a square then $\Sigma$ -
 generic
bundles exist. Since being $\Sigma$ - generic is an open condition,
this will imply that such bundles form an open dense subset in the
moduli space.
We shall also see  that for (some other) open dense subset in the
moduli space the variety $\Sigma$ is indeed a surface. However,
there are important particular cases when $\Sigma$ is reducible
and contains components of higher dimension, see n. 3.5 below.

\vskip .3cm

\noindent {\bf 2.6.} Consider now the partial transpose $a_V$ of
the tensor $a$ given in (2.4). We define the determinantal variety
$Z\i P(V^*)$ by
$$Z = \{\lambda \in P(V^*): \quad {\rm rank}\,\,\, a_V(\lambda) \leq n-2\}.
\eqno (2.7)$$
It will be important for us to consider $Z$ as a scheme with the scheme
structure given naturally by (2.7). This means that we choose bases in $H$
 and $M$ and regard $a_V$ as a $(n-1)\times n$ -matrix whose entries are
linear forms in $\lambda$. The $n$ maximal minors of this matrix are
taken to be the equations of the subscheme $Z$.

\vskip .2cm

Since Hom$(H,M)_{n-2}$ has codimension 2 in Hom$(H,M)$, generically one
expects $Z$ to be 0-dimensional
and reduced. If this is indeed the case, we shall call
the tensor $a$ (and the bundle $E$) $Z$ - {\it generic}.
 It follows from [Hu 1]
that $Z$ -generic bundles exist for any values of $n$.
Namely, the so-called
Hulsbergen bundles will be $Z$ -generic (see also \S 4
 for discussion of
these bundles). Thus $Z$ -generic bundles form an open
 dense subset
in the moduli space.

If $a$ is $Z$-generic then, by Proposition 2.4. c), the degree of the 0-
dimensional scheme $Z$ equals  deg Hom$(H,M)_{n-2} = {n\choose 2}$.
Moreover, the multiplicity of any point $\lambda\in Z$ in $Z$ is at least
$n - r(\lambda) \choose 2$ where $r(\lambda) =
{\rm rank}\,\,\, a_V(\lambda)$

The meaning of $Z$ is as follows.

\proclaim 2.7.~Lemma. The support of the scheme $Z$ is precisely the set
of jumping lines of $E$.

\noindent {\sl Proof:} This immediately follows from
considerations of n.1.3.

\vskip .3cm

\noindent {\bf 2.8.} Let ${\cal J}_Z \i {\cal O}_{P(V^*)}$ be the
sheaf of
ideals of the subscheme $Z$. By construction of $Z$ (see n. 2.6),
maximal
minors of the $(n-1) \times n$ - matrix $a_V$ are global sections of
${\cal J}_Z(n-1)$. In invariant terms, we consider the
linear map
$$a_V^*: H\otimes M^* \rightarrow V \eqno (2.9)$$
and, by taking its $(n-1)$ -st symmetric power, we get a
linear map
$$\bigwedge^{n-1}H \otimes \bigwedge^{n-1}M^* \hookrightarrow
S^{n-1}(H\otimes M^*) \longrightarrow S^{n-1}V =
 H^0\bigl (P(V^*), {\cal O}(n-1)\bigl)
\eqno (2.10)$$
whose image is contained in $H^0\bigl(P(V^*),
{\cal J}_Z(n-1)\bigl)$.
 It will be convenient for
us to rewrite (2.10) as
$$A: M\otimes \biggl(\bigwedge^{n-1}H \otimes \bigwedge^{n}M^*\biggl)
\quad \longrightarrow
\quad H^0\bigl(P(V^*), {\cal J}_Z(n-1)\bigl).\eqno (2.11)$$
The 1-dimensional vector space $\bigwedge^{n-1}H \otimes
 \bigwedge^{n}M^*$
can be chased away by choosing bases in $H$ and $M$.

\proclaim 2.9.~Proposition. If {\rm dim} $Z = 0$ then the
operator $A$ in (2.11)
is an isomorphism. In other words, the linear system
of curves of degree $n-1$
through $Z$ is generated by maximal minors of $a_V$.

\noindent {\sl Proof:} We associate to $a_V$,
in a standard way, a morphism
$\tilde a$ of sheaves on $P(V^*)$ and denote its
cokernel by ${\cal F}$:
$$0\rightarrow H\otimes {\cal O}_{P(V^*)}(-1)
\buildrel \tilde a
\over\longrightarrow M\otimes {\cal O}_{P(V^*)} \rightarrow {\cal F}
\rightarrow 0 \eqno (2.12)$$
(the fact that $\tilde a$ is injective, follows from dim $Z = 0$). We
claim the following:

\proclaim 2.10.~Lemma. ${\cal F}$ is isomorphic to ${\cal J}_Z(n-1)$.
Under this isomorphism the natural map $M \rightarrow H^0
\bigl({\cal J}_Z(n-1)\bigl)$
corresponds, up to a scalar multiple, to the map $A$ from (2.11).

Clearly, Lemma 2.10 implies our proposition in virtue of the exact
cohomological sequence of (2.12).

\vskip .2cm

\noindent {\sl Proof of the lemma:}
We choose a bases $h_1,...,h_{n-1}\in H$
and $m_1,...,m_n \in M$. This makes it possible to speak about
 the determinant
$\det [v_1,...,v_n]$ of a system of $n$ vectors in $M$ (this
is just $|b_{ij}|$
where $v_i = \sum b_{ij}m_j$). We define a morphism of sheaves
$\psi: {\cal F} \rightarrow {\cal J}_Z(n-1)$ i.e. a morphism $\Psi:
M\otimes {\cal O}_{P(V^*)}\rightarrow {\cal J}_Z (n-1)$ vanishing on
Im$(\tilde a)$, as follows. Let $m = m(\lambda)$ be a local section of
$M\otimes {\cal O}_{P(V^*)}$ i.e.      an $M$ - valued function in
$\lambda$ homogeneous of degree 0. We put $\Psi(m)$ to be the
homogeneous (of degree $n-1$) function
$$\lambda \mapsto \det\bigl[ m(\lambda), a_V(\lambda)(h_1), ...,
a_V(\lambda)(h_{n-1})\bigl].$$
This defines $\psi$. It is clear that $\psi$ is injective. The fact that
$\psi$ is surjective follows by comparing Chern classes of ${\cal F}$
and ${\cal J}_Z(n-1)$. The rest of the lemma is obvious.

\vskip .3cm

\noindent {\bf 2.11.} We continue to assume
 that dim $Z = 0$. Let $S$ be the
blow up of $P(V^*)$ along $Z$ and $\pi_S: S \rightarrow P(V^*)$ be the
canonical projection. In virtue of Proposition 2.9 the linear system
of curves of degree $n-1$ through $Z$ defines a regular map
$p:S\rightarrow P(M^*)$. A generic point $s = \pi_S^{-1}(\lambda)
\in S,\,\,
\lambda \in P(V^*)$ goes under $p$ into the hyperplane in $P(M)$
consisting of $m$ such that $A(m) \in H^0\bigl(P(V^*),
{\cal J}_Z(n-1)\bigl)$
vanishes at $\lambda$ as well. Here $A$ is as in (2.11).
The interpretation
of $A$ in Lemma 2.10 shows that $p(S)$ is contained in the
determinantal
variety $\Sigma \i P(M^*)$ as an irreducible component.
We shall denote variety $p(S)$ (typically a surface) by
$\Sigma' \i \Sigma$.

\vskip .3cm

\noindent {\bf 2.12.} Suppose that our bundle $E$ is $\Sigma$ -
generic. Then we have a regular map
$$\pi_\Sigma: \Sigma \rightarrow P(V^*)$$
which takes $\mu\in\Sigma\i P(M^*)$ to the linear
 subspace Im$(a_M(\mu))\i V$
(this subspace has dimension 2 by the assumption
of $\Sigma$ -genericity). The
map $\pi_\Sigma$ is the analog of the blow-down of
 a cubic surface onto a plane.

If the bundle $E$ is not $\Sigma$ -generic, the map
$\pi_\Sigma$ will be defined
on the open part $\Sigma_0\i \Sigma$ consisting of
 $\mu$ such that
$a_M(\mu)$ has rank exactly 2.

\vskip .2cm

For $\lambda\in P(V^*)$ the fiber $\pi_\Sigma^{-1}(\lambda)$
is the projective
space $P({\rm Ker} \,\,\, a_V(\lambda)^*)$. The dimension of
this fiber is equal to $n - {\rm rank} \,\,\, a_V(\lambda) - 1$.
 Hence
$\pi_\Sigma$ is an isomorphism over the complement of Supp$(Z)$.
On the other hand,
if the rank of $a_V(\lambda)$ is small the fiber
$\pi_\Sigma^{-1}(\lambda)$ will
have dimension $\geq 2$ and the variety $\Sigma$ will
be reducible. We shall
see in \S 3 that such situations do occur for stable bundles.

\vskip .3cm

\proclaim 2.13.~Proposition. Assume that $E$ is $Z$ - generic and
no $n-1$ points of $Z$ lie on a line. Then:
\item{(a)} The map $p: S \rightarrow \Sigma$ is an isomorphism
(so, in particular, $\Sigma' = \Sigma$);
\item{(b)} $\Sigma$ is a projectively Cohen - Macaulay
surface in $P(L)$ of degree $(n-1)^2 - {n\choose 2}$.

\noindent {\sl Proof:} By the exact sequence (2.12), we have
$$h^0({\cal O}_{P(V^*)}(n-1)) - h^0({\cal J}_Z(n-1)) =
(1/2) n(n+1) - n = (1/2)n(n-1) = $$
$$ = h^0({\cal O}_{P(V^*)}(n-2))
- h^0({\cal J}_Z(n-2)) >$$
$$> h^0({\cal O}_{P(V^*)}(n-3)) - h^0({\cal J}_Z(n-3)) =$$
$$ =
(1/2)(n-1)(n-2).$$
Denoting by $H(Z,n) = h^0({\cal O}_Z(n))$ the Hilbert function
of the subscheme $Z$, we obtain that
$$n-1 = \min \{ t: H(Z,t) = H(Z,t-1) \}.$$
By [DG], this implies that the linear system of curves of degree
$n-1$ through $Z$ maps $S = {\rm Bl}_Z(P^2)$ isomorphically
into $P(H^0({\cal J}_Z(n))^*) = P(M^*)$. By [Gi] the image of
this map i.e., the variety $\Sigma'$, is
projectively Cohen - Macaulay. Recall that this means that
the projective coordinate ring of $\Sigma'$ is Cohen - Macaulay.
In particular, we get that $\Sigma'$ is projectively normal
i.e., for any $k\geq 0$ the restriction map
$$H^0\bigl(P(M^*), {\cal O}(k)\bigl)
 \longrightarrow H^0\bigl(\Sigma', {\cal O}(k)\bigl)$$
is surjective. Since the rational map $P(V^*) \rightarrow
\Sigma$ is given by the linear system
of curves of degree
$n-1$ through $Z$, we obtain the assertion about the degree of $\Sigma'$.
Since $E$ is $Z$ - generic, the fiber of the map
$\pi_\Sigma: \Sigma \rightarrow P(V^*)$ over each point $z\in Z$
is isomorphic to $P^1$. Since $\Sigma'$ and $\Sigma$ coincide
outside the union of the fibers $\pi_\Sigma^{-1}(z)$,
$z\in Z$, this implies that $\Sigma' = \Sigma$. Q.E.D.

\vskip .3cm

\noindent {\bf 2.14.}
Let $z\in Z \i P(V^*)$.
We denote the fiber
$$\pi_\Sigma^{-1}(z) = P \bigl({\rm Ker}  (a_V(z)^*\bigl)
\i P(M^*) \quad {\rm by} \quad A_z.$$
The corresponding linear subspace ${\rm Ker}  (a_V(z)^*) \i M^*$ of
which $A_z$ is the projectivization, will be denoted by ${\bf A}_z$.

Consider the space
$$H_z = {\rm Ker} \,\,\, a_V(z) \i H.$$
We also consider the linear subspace
$$ {\bf A}'_z = \bigcap_{h\in H_z} {\rm Ker} \,\,\, a_H(h)
\quad \i \quad M^*$$
and denote its projectivization by $A'_z \i P(M^*)$.

\vskip .2cm

The collection of projective subspaces $A_z, A'_z, \,\, z\in  Z$,
 forms the analog
of a Schl\"afli double - six on a cubic surface in $P^3$.

\vskip .2cm

In our case  $A'_z$
 lies
on $\Sigma$ but $A_z$ does not, in general, do so. Indeed,
the typical situation (see Proposition 2.13) is that $\Sigma$
is  a surface, that for any $z\in Z$
we have rk$(a(z)) = n-2$ and so dim $A_z = 1$, dim $A'_z = n-3$.
 So for $n\geq 5$ $A'_z$ cannot lie on $\Sigma$. The
relation of $A'_z$ with the component $\Sigma' = p(S) \i \Sigma$
is as follows.

\proclaim 2.15.~Proposition. Assume that $E$ is $Z$ - generic.
Then $A'_z$ is a subspace of codimension 2 in $P(M^*)$
which intersects the surface $\Sigma'$ along a curve. The image
of this curve under the projection $\pi_\Sigma:
\Sigma \rightarrow P(V^*)$ is the unique curve of degree $n-2$
which passes through the points $z', z\in Z$

For the case $n=4$ we get the standard description of the
second sextuple of lines of the double -  six on the
cubic surface
as the inverse image of quadrics containing some 5 of the 6 points
of $Z$. For $n \geq 5$ instead of the property that $A'_z$ lies
 on $\Sigma'$ we have that $A'_z \cap\,\, \Sigma'$
is a curve (instead of a set of isolated points, as one would
expect by dimension count).

\vskip .2cm

\noindent {\sl Proof:} Since the rank of $a_V(z)$ equals $n-2$,
we have dim$(H_z) =1$. Thus ${\bf A}'_z$ is the kernel of the
map $a_V(h): M^* \rightarrow V$ where $h$ is any non-zero
vector from $H_z$. Note that the rank of this map equals 2.
In fact, otherwise $Z$ would contain a line as an irreducible
component. This shows that dim$({\bf A}_z) = n-2$. Now
let us observe that $A'_z = P({\bf A}'_z)$ intersects
each $A_{z'}$ for $z'\neq z$. Indeed, the sum of linear subspaces
${\bf A}'_z + {\bf A}_{z'}$ is contained in the hyperplane of
zeroes of the linear form $a(h,z') \in M = (M^*)^*$, where $a$
is as in (2.1).

Let
$\{H(\lambda)\}_{\lambda\in P^1}$
be the pencil of hyperplanes in $P(M^*)$ which
contain the subspace $A'_z$. It cuts out a pencil ${\cal P}$
of curves on $\Sigma$ with the base locus $A'_z\cap\Sigma$.
For each $z'\neq z$ one of the hyperplanes $H(\lambda)$
contains the line $A_{z'}$. Thus each $A_{z'}$ contains
one of the base points of the pencil ${\cal P}$.
Under the rational map $P(V^*) \rightarrow P(M^*)$
(given by curves of degree $n-1$ through $Z$) the preimage of the
pencil $\{H(\lambda)\}$ is some pencil of curves of degree
$n-1$ passing through $Z$. Let ${\cal C}$ be its moving part
and $F$ be its fixed curve. Let $d$ be the degree of $F$
(zero if $F = \emptyset$). Curves of the pencil ${\cal C}$
have degree $n-1-d$.
Suppose that they pass through some $m$ points
say, $z_1,...,z_m$ of $Z$. Then,
since $z_i$ remain basic for ${\cal C}$ after the blow - up,
 curves from ${\cal C}$ have the same tangent
direction at each $z_i \neq z$.
 The curve $F$ passes through the remaining
$(1/2)n(n-1) - m$ points of $Z$. Consider a typical curve $C \in
{\cal C}$. Let $\tilde C$ be its proper transform in $S =
{\rm Bl}_Z(P(V^*))$. Since $\tilde C$ moves, its self -
intersection index is non - negative so we get
$$0 \leq \tilde C^2 \leq (n-d-1)^2 - 2(m-1) - 1 = (n-d)(n-d-1) -
2m - (n-d-2).$$
If $n-d-2 \geq 0$, we obtain that $(n-d)(n-d-1) - 2m \geq 0$
thus there exists a plane curve of degree $n-d-2$ passing through
$z_1,...,z_m$. Together with the curve $F$, it defines
a curve of degree $n-2$ passing through all the points of $Z$.
But Lemma 2.10 and the exact sequence (2.12) show
 that this is impossible.
So we must have $d = n-2$ and hence $m=1$, so $F$ is a curve
of degree $n-2$ which passes through all the points of $Z$
except $z$. If there is another curve, say, $F'$, with this
property
then we would have a pencil of curves of degree $n-2$
through $Z - \{z\}$. This pencil must then contain a curve
passing also through $z$. This, as we have just seen,
is impossible.

\vskip .3cm

\noindent {\bf 2.17.} Up until now we worked exclusively with
the tensor $a$ from (2.1). Now we take into account the
 non-degenerate quadratic form $B\in S^2M^*$
from (2.2). Let $C = B^{-1}$ be the inverse
quadratic form
on $M^*$. The following result justifies the name "Schur quadric"
 for
the quadric defined by $C$.

\proclaim 2.17.~Theorem. Let $z\in {\rm Supp} (Z)$. Then $A_z$ is contained
in the orthogonal complement $\bigl(A'_z\bigl)^\bot_C$
of $A'_z$
 with respect to  $C$ . If,
moreover, rk $a_V(z) = n-2$ then we have equality $A'_z =
\bigl(A_z\bigl)_C^\bot$.

\noindent {\sl Proof:} For any $\lambda\in V^*$ let
 $$b(\lambda) = a(\lambda)^* \circ B:\quad M \longrightarrow H^*$$
where $a(\lambda)$ is the map induced by the $a$ from (2.1).
Then
$$B^{-1}({\bf A}'_z) =
\bigl\{m\in M: (b(\lambda)(m), h) = 0,\quad \forall
\lambda\in V^*, h\in H_z = {\rm Ker} \,\,\, (a_V(z))\bigl\}.$$
For any $m\in {\bf A}_z^\bot = a_V(z)(H)$ we write $m = a(z)(h')$ for
some $h'\in H$ and obtain
$$(b(\lambda)(a_V(z)(h')), h)
\quad =  \quad (b(\lambda)(a_V(z)(h), h') = (0,h') = 0.$$
Here we use the property $(\alpha 3)$ from n.1.1.
 Thus we obtain
$${\bf A}_z^\bot \i B^{-1}({\bf A}'_z).$$
If rank$\,(a(z)) = n-2$ then dim ${\bf A}_z = 2$,
dim $H_z = 1$ and dim ${\bf A}'_z = n-2$.
Thus the dimensions of the spaces ${\bf A}_z^\bot$ and
 $B^{-1}({\bf A}'_z)$ are the same
so these spaces are equal.  Theorem is proven.

\vskip .3cm

\noindent {\bf 2.18. Remark.} Let $Z$ be any set of $n\choose 2$
points in $P^2$ such that no curve of degree $n-2$
contains $Z$ and no lines pass through $n-1$ points of $Z$.
The linear system of curves of degree $n-1$ through $Z$
defines a rational map of $P^2$ into $P^{n-1}$
whose image is a nonsingular surface $X$ classically known as a
{\it White surface} [R]. If $n=4$, this is a cubic surface.
The surface $X$ is given by vanishing of maximal minors
of a $3\times (n-1)$ matrix of linear forms. A modern proof of
these results can be found in [DG] and [Gi].

\vskip .1cm

Every White surface comes equipped with a set of $n\choose 2$ lines
$E_z, z\in Z$ corresponding to exceptional curves of the blow - up
Bl$_Z(P^2)$ and a set of $n\choose 2$ curves $C_z$ of degree
$(n-2)(n-4)/2 + 1$. The curve $C_z$ is the image  of the
(unique)  plane curve
of degree $n-2$ passing through $Z - \{z\}$. Each curve $C_z$
spans a subspace $E'_z$ of codimension 2 in $P^{n-1}$.
We have $E_z \cap \,\, E'_z = \emptyset$ but $E_z\cap \,\, E'_{z'}
\neq \emptyset$ for $z'\neq z$. This situation is analogous to
a configuration of a double - six on a cubic surface.

\vskip .1cm

Propositions 2.13 and 2.15 imply that for a $\Sigma$ - generic
stable bundle $E$ the variety $\Sigma$ is a White
surface. However, by counting constants it follows
that not every White surface comes in this way, as soon as $n \geq 5$.
Although one can reconstruct a linear map
$$a: H \otimes V^* \cong {\bf C}^{n-1} \otimes {\bf C}^3 \,\,\,
\longrightarrow  \,\,\, M = {\bf C}^n$$
 from a
 determinantal representation of $X$, there does not exist, in
general, a quadratic form $B$ on $M$ such that $a$ satisfies
the property $(\alpha 3)$ from n.1.1. By Theorem 1.2 the existence
of such a $B$ is necessary and sufficient in order that $X = \Sigma$
for some $Z$ - generic bundle $E$. It seems likely that these
conditions are equivalent to the existence of a "Schur quadric"
for the "double - six" $\{E_z, E'_z\}$ i.e., a quadric $Q$ in $P^{n-1}$
such that $E_z$ and $E'_z$ are orthogonal with respect to the
(quadratic form defining) $Q$.

\vskip .3cm

\noindent {\bf 2.20.} The role of the Schur quadric $Q$ (see n.2.2)
in the description of
jumping lines of the second kind is given by the following
remark.

\proclaim 2.21.~Proposition. Let $\Sigma_0 \i \Sigma$ be the open set of
$\mu$ such that the rank of $a_V(\mu)$ equals 2 (so $\Sigma_0 = \Sigma$
if the bundle is $\Sigma$ - generic). Let $\pi_\Sigma: \Sigma_0 \rightarrow
P(V^*)$ be the projection defined in n. 2.12. Then the curve $C(E)$
of jumping lines of second kind coincides with the closure
of $\pi_\Sigma (Q\cap \Sigma_0)$.

In particular, when the bundle $E$ is $\Sigma$ - generic, we have
$C(E) = \pi_\Sigma (Q\cap\Sigma)$ .

\noindent {\sl Proof:} This is a reformulation of what has been done
in n.1.5.

\vskip .3cm

As an application of our formalism  of Schur quadrics let us
prove a statement about the singular tangent lines of the curves
of JLSK which strengthens, under assumptions of genericity,
a theorem of Hulek. More precisely, Hulek [Hu 1] has proven the
following fact.

\proclaim 2.22.~Theorem. Let $l\in C(E)$ be a JLSK of $E$.
Suppose that $E|_l \cong {\cal O}(-1-k) \oplus {\cal O}(k)$
with $k\geq 1$.
Then $l$ is a singular point of the curve $C(E)$ of multiplicity
$2k$ and for any line $T$ in the tangent cone of $C(E)$
at $l$ the intersection index of $C(E)$ and $T$
at $l$ is at least $2k+2$.

Assume that $E$ is $Z$ - generic. Then every singular
point $l$ of $C(E)$ is a double point
(a node or, possibly, a  cusp of type $y^2 = x^r$).
Theorem 2.22 gives
that in this case there exist at least one line
$T$ through the point $l$ with intersection index $\geq 4$.

\vskip .2cm

We claim that the case of the cusp does not occur for $Z$-generic $E$.
Call an ordinary
double point $p$ of a plane curve $C$ a {\it biflexnode} if each
of the two branches
has a flex at this point i.e. each of the two tangents has
the intersection index $\geq 4$ with $C$ at $p$.

\proclaim 2.23.~Theorem. Assume that the bundle $E$ is $Z$ -
generic. Then every singular point of $C(E)$ is a biflexnode.

\noindent {\sl Proof:} Let $z$ be a singular point of $C(E)$.
Then $z\in Z$. The branches of $C(E)$ at $z$ correspond to
the points of intersection of the line $A_z$ and the Schur
quadric $Q$. Note that $Q$ cannot be tangent to $A_z$
since otherwise we would have $A_z \cap A'_z \neq \emptyset$
which contradicts Proposition 2.15.
 This proves that the point $z$ is an ordinary
node.

Although Theorem 2.22 allows us to finish the proof, we prefer to
give an independent proof based on the properties of the Schur quadric.

Now let $x$ be one of the two points of
$Q\cap A_z$ and let $\Pi$ be a hyperplane in $P(M^*)$
which is spanned by the point $x$ and the codimension 2
subspace $A'_z$. Let $Q(z)$ denote the quadric in $A'_z$ cut out
by $Q$. For any point $y\in Q(z)$ the line $<x,y>$ is contained
in $Q$. This implies that $\Pi$ is tangent to $Q$ at $x$.
Let $\tilde C (E) = \Sigma \cap \,\,Q$ be the proper inverse
transform of the curve $C(E)$ in $\Sigma$, under the blow-down
$\pi_\Sigma: \Sigma \rightarrow P(V^*)$. Let $\tau$
be the tangent line to $C(E)$ at $z$ at the branch corresponding
to $x$ and let $\tilde \tau$ be its proper inverse transform on
$\Sigma$.

 Under the correspondence between hyperplanes in
$P(M^*)$ and curves of degree $n-1$ in $P(V^*)$ through $Z$,
the hyperplane $\Pi$ corresponds to the reducible curve
$\tau + C_z$ where $C_z$ is the plane curve of degree $n-2$
passing through $Z - \{z\}$.
This implies that $\tilde \tau \i \Pi$ and so
$$T_x (\tilde\tau) = \Pi \cap \,\, T_x(\Sigma) =
T_x(Q) \cap \,\, T_x (\Sigma) = T_x (\tilde C(E)).$$
This shows that $\tilde\tau$ is tangent to $\tilde C(E)$ at the
point $x$. Obviously this implies that $\tau$ is a flex
tangent at the branch of $C(E)$ at $z$ corresponding to $x$.
Theorem is proven.

\beginsection

\centerline {\bf \S 3. Logarithmic bundles.}

\vskip 1cm

\noindent {\bf 3.1.} Consider a projective plane $P^2 = P(V)$, dim $V = 3$.
Let ${\cal H} = (H_1,...,H_m)$ be an arrangement of $m$ lines
in $P(V)$ in general position (i.e., no three of these lines have a common
point). Let $E({\cal H}) = \Omega^1_{P(V)}(\log {\cal H})$ be the sheaf of
1-forms on $P(V)$ with logarithmic poles along $H_i$. Since ${\cal H}$
is a divisor with normal crossings, $E({\cal H})$ is locally free i.e.
we can and will regard it as a rank 2 vector bundle. It was proven in [DK]
that this bundle is stable.

We further suppose that the number of lines is even: $m=2d$.
 In this case $c_1 E({\cal H}) = 2d-3$. The normalized bundle
 $E_{\rm norm}({\cal H}) = E({\cal H})(-d+1)$
is a stable bundle with $c_1 = -1, c_2 = (d-1)^2$. In this section we apply
considerations of \S\S 1,2 to bundles $E_{\rm norm}({\cal H})$.

\vskip .3cm

\noindent {\bf 3.2.} It was shown in [DK] that the bundle $E({\cal H})$
has a resolution of the form
$$0\rightarrow I\otimes {\cal O}_{P(V)}(-1) \buildrel \tau \over
\longrightarrow W\otimes {\cal O}_{P(V)} \rightarrow E({\cal H})
\rightarrow 0.\eqno (3.1)$$
In (3.1) the space $W$ is defined as
$$W = \biggl\{ (a_1,...,a_{2d}) \in {\bf C}^{2d}: \sum a_i = 0
\biggl\}.\eqno
(3.2)$$
The space $I$ is defined as follows. Let $f_i \in V^*$ be a linear
equation of the line $H_i$. Then $I$ is the space of relations
among $(f_1,...,f_{2d})$ i.e.,
$$I = \biggl\{(a_1,...,a_{2d}) \in {\bf C}^{2d}: \sum a_i f_i = 0
\biggl\}. \eqno (3.3)$$
The map $\tau$ is induced by the canonical map
$$t: I\otimes V \rightarrow W, \quad
(a_1,...,a_{2d}) \otimes v \mapsto \bigl(a_1f_1(v),...,
a_{2d}f_{2d}(v)\bigl) \eqno
(3.4)$$
called the fundamental tensor of ${\cal H}$.

\vskip .3cm

\noindent {\bf 3.3.} By twisting the resolution (3.1) with ${\cal O}(-d+1)$
we get a resolution for $E_{\rm norm}({\cal H}) = E({\cal H}) (-d+1)$. From
this
it is immediate to find the data defining the Hulek's monad for
$E_{\rm norm}({\cal H})$ (see \S 1). To formulate the answer neatly, let again
$f_j\in V^*$ be the equation of $H_j$. For any $m\geq 1$
denote by $\partial/\partial f_j:
S^mV \rightarrow S^{m-1}V$ the derivation corresponding to $f_j$ regarded
as a constant vector field on $V^*$.  We define the following map
$$t_{(m)}: S^mV \otimes I \rightarrow S^{m-1}V \otimes W,\eqno (3.5)$$
$$p\otimes (a_1,...,a_{2d}) \mapsto
\biggl(a_1 {\partial p\over\partial f_i},
.., a_{2d} {\partial p\over\partial f_{2d}}\biggl)\eqno (3.6)$$
where we regard $S^{m-1}V \otimes W$ as the space of collections
$(q_1,...,q_{2d})$ of polynomials $q_j\in S^{m-1}V$ summing up to 0.

Now the vector spaces in the monad for
 $E_{\rm norm}({\cal H})$ have the form
$$H = H^1 \bigl(E_{\rm norm}({\cal H})(-2)\bigl) =
 H^1\bigl(E({\cal H})(-d-1)\bigl) = {\rm Ker} \,\,
(t_{(d-1)}); \eqno (3.7)$$
$$M = H^1 \bigl(E({\cal H})(-d)\bigl)
= {\rm Ker}\,\,(t_{(d-2)});\eqno (3.8)$$
$$H' = H^1 \bigl(E({\cal H})(-d+1)\bigl)
 = {\rm Ker}\,\, (t_{(d-3)}),\eqno (3.9)$$
as it follows immediately  from  the
resolution (3.1). For example, the map
$t_{(d-1)}: S^{d-1}V \otimes I \rightarrow
S^{d-2}V \otimes W$ in (3.7) appears
as the map
$$H^2 \bigl(P(V), {\cal O}(-d-2) \otimes I\bigl)
\rightarrow H^2 \bigl(P(V), {\cal O}(-d-1)
\otimes W\bigl)$$
in the long exact sequence of cohomology of the resolution (3.1)
tensored with ${\cal O}(-d-1)$.

As regards maps in the monad (1.1), we shall only need the explicit form of
the operator
$$b_M: M\rightarrow V\otimes H' \eqno (3.10)$$
defined by the map $b$ in (1.1). Namely, $b_M$ is induced by
$$\psi\otimes {\rm Id}_I: S^{d-2}V \otimes I \rightarrow V\otimes S^{d-3}V
\otimes I \eqno (3.11)$$
where $\psi: S^{d-2}V \rightarrow V\otimes S^{d-3}V$ is the canonical
$GL(V)$ - equivariant
embedding.
The map $a$ in (1.1) is dual to $b$ by means of the form $B$.

The following is the main result of this section.

\proclaim 3.4.~Theorem. Any bundle $E_{\rm norm}({\cal H})$ is $\Sigma$
-generic
(see n. 2.5).

\noindent {\sl Proof:} In the notation of \S 2 we have to prove that
$$a_M(M^*)\quad  \cap \quad {\rm Hom}(V^*, H^*)_1
\quad\quad = \quad\quad \{0\}.\eqno (3.12)$$
We have a commutative diagram
$$\matrix{M^* & \buildrel a_M \over\longrightarrow & {\rm Hom}(V^*, H^*)\cr
\big\downarrow&&\big\downarrow\cr
M&\buildrel b_M\over\longrightarrow & {\rm Hom}
(V^*, H')}$$
where the left vertical arrow is induced by the form $B$
and the right vertical arrow --- by the isomorphism $H^* = H'$ (see n. 1.2).
It is enough therefore to prove that
$$b_M(M) \cap {\rm Hom}(V^*, H')_1 = \{0\}. \eqno (3.13)$$
Let $m = \sum p_i \otimes x_i$ be an element of $M \i S^{d-2}V \otimes I$,
so $p_i \in S^{d-2}V, x_i \in I$. The element $m$ is mapped by $b_M$ into an
element of ${\rm Hom}(V^*, H')_1$ if and only if there is $v\in V$ such that
each $p_i$ equals $vq_i$ for some $q_i \in S^{d-3}V$ and also
$\sum q_i\otimes x_i \in H'$. Each $x_i\in I$ is in fact a vector
$x_i = (x_i^{(1)}, ..., x_i^{(2d)})$ such that $\sum_{j=1}^{2d} x_i^{(j)}
f_j = 0$. Since $m$ belongs to $M = {\rm Ker} \,(t_{(d-1)})$, we have, by (3.5)
and (3.6):
$$\sum_i x_i^{(j)} {\partial (vq_i)\over\partial f_j} = 0,\quad j=1,...,2d.
\eqno (3.14)$$
By applying Leibnitz' rule for $\partial/\partial f_j$ and taking into account
the fact that $\sum_i q_i\otimes x_i \in H' = {\rm Ker}(t_{(d-3)})$, we get
the equalities
$$f_j(v)\sum_i x_i^{(j)} q_i = 0,\quad j=1,...,2d.\eqno (3.15)$$
We claim that these equalities imply that $q_i=0$ for all $i$. Indeed, let
$\lambda: S^{d-3}V \rightarrow {\bf C}$ be any linear functional. Consider
the vector
$$y = \sum_i \lambda(q_i)x_i \in I.$$
If we write $y$ in terms of its components: $y = (y^{(1)}, ..., y^{(2d)})$
then (3.15) implies that
$$f_j(v)y^{(j)} = 0,\quad j=1,...,2d.$$
Let $J = \{j: f_j(v) = 0\}$. Since the lines $\{f_j=0\}$ are in general
position, $|J|\leq 2$. For $j\notin J$  we have therefore $y^{(j)} = 0$.
Since $y\in I$, we have
$$ 0 = \sum_{j=1}^{2d} y^{(j)}f_j = \sum_{j\in J} y^{(j)}f_j$$
which means that we have a nontrivial linear relation among $|J|\leq 2$
elements of \hfill\break $\{f_1,...,  f_{2d}\}$. This contradicts the general
position
of $\{f_i=0\}$ so the vector $y\in I$ is zero. In other words, for any
linear functional $\lambda: S^{d-3}V \rightarrow {\bf C}$ we have
$\sum \lambda(q_i) x_i = 0$ in $I$. This means that $\sum q_i \otimes x_i
= 0$ in $S^{d-3}V \otimes I$ and Theorem 3.4 is proven.

\vskip .3cm

\noindent {\bf 3.5.} Let $Z$ be the subscheme of jumping lines of
$E_{\rm norm}({\cal H})$. As was shown in [DK] (Proposition 7.4),
 the lines $H_i$ belong to
$Z$. Moreover,
$$E_{\rm norm}({\cal H})|_{H_i} = {\cal O}_{H_i}(1-d) \oplus {\cal O}_{H_i}
(d-2).\eqno (3.16)$$
Denote, as usual, by $f_i \in V^*$ the equation of $H_i$.
The equality
(3.16) means that the matrix $a_V(f_i)$ (see formula (2.4))
has rank $n-d-1$.
By Proposition 2.4 d) this implies that the multiplicity of
each $H_i$ as a point of $Z$ is at least $(d-1)(d-2)/2$.
The total degree of $Z$, however,
equals to $n\choose 2$ where $n = c_2(E_{\rm norm}({\cal H})) = (d-1)^2$.
Thus one expects that for $d\geq 4$ there will be many other jumping lines
apart from $H_1,....,H_{2d}$.

Let us also note that the fibers of the map $\pi_\Sigma:
\Sigma \rightarrow P(V^*)$
introduced in n. 2.12 over points $H_i \in P(V^*)$ are projective spaces
of dimension $d-2$. This means that for $d\geq 4$ the
determinantal  variety $\Sigma$
("cubic surface") will be always reducible.

\beginsection

\centerline {\bf \S 4. Examples.}

\vskip 1cm

\noindent {\bf 4.1.} In this section we shall illustrate geometric
constructions of \S 2 on some particular classes of bundles. The example
with cubic surfaces and Schur quadrics (which motivated the present paper)
will be considered in n.4.4.

In each of the examples below we shall indicate the value of $n = c_2$
(we assume $c_1 = -1$) and describe the following geometric objects
(all introduced in \S 2):

\vskip .2cm

\item{a)} The subscheme $Z\i P(V^*)$ of jumping lines. If dim $Z = 0$ then deg
$Z =
{n\choose 2}$.

\vskip .1cm

\item{b)} The determinantal variety $\Sigma \i P(M^*)$ (the analog of the cubic
surface). It comes with a natural map $p: {\rm Bl}_ZP(V^*) \rightarrow \Sigma$
whose
image is a component of $\Sigma$. The map $p$ is given by the linear system
of curves of degree $n-1$  in $P(V^*)$ through $Z$.

\vskip .1cm

\item{c)} The Schur quadric $Q\i P(M^*)$.

\vskip .1cm

\item{d)} The curve $C(E)\i P(V^*)$ of JLSK. Its degree is $2n-2$. It can be
described
as $\overline {\pi_\Sigma(\Sigma_0\cap Q)}$ where $\pi_\Sigma:
 \Sigma_0 \rightarrow P(V^*)$
is the projection of the generic part of $\Sigma$ introduced in n. 2.12.

\vskip .1cm

\item{e)} The
projective subspaces $A_z, A'_z, z\in Z$
(the analog of the double - six).

\vskip .2cm

\noindent By $M(-1,n)$ we shall denote the moduli space of stable
rank 2 vector bundles on $P^2$ with $c_1 = -1, \,\, c_2 = n$.
It is an irreducible variety of dimension $4n-4$, see [Hu 1] [OSS].

\vskip .3cm

\noindent {\bf 4.2. The case $n=2$.} This case was considered by Hulek [Hu 1].
The features are as follows:

a) $Z$ consists of just one point $z_0\in P(V^*)$. This point corresponds to
the 1-dimensional kernel of
$$a_V: V^* \rightarrow {\rm Hom}(H,M) = {\bf C}^2.$$

b) The determinantal variety $\Sigma \i P(M^*) = P^1$ coincides with $P(M^*)$.
The regular map $p: {\rm Bl}_ZP(V^*) \rightarrow
 \Sigma$ is the natural projection
${\rm Bl}_{z_0}P^2 \rightarrow P^1$.

c) The Schur quadric $Q\i P(M^*) = P^1$ consists of two distinct points.

d) The curve $C(E)$ is $\pi(p^{-1}(Q))$ where $\pi:
{\rm Bl}_ZP(V^*) \rightarrow
P(V^*)$ is the projection. In other words, $C(E)$ is the union of two
distinct lines through $z_0$.

e) The "double - six" is as follows: $A_{z_0} = P(M^*), A'_{z_0} =
\emptyset$.

\vskip .3cm

\noindent {\bf 4.3. The case $n=3$.} There may be several possibilities for
$Z$ which were also
listed by Hulek [Hu 1]. We shall consider only the most
 generic case when $Z$ consists of three
distinct non-collinear points. In this case the features are as follows:

\vskip .2cm

b)The variety $\Sigma\i P(M^*)$ again coincides with $P(M^*) = P^2$.
The regular map ${\rm Bl}_Z P(V^*) \rightarrow P(M^*) = \Sigma$
resolves  the standard
Cremona transformation $c: P(V) = P^2 \rightarrow P^2 = P(M^*)$ defined
by quadrics through $Z$ (three points).
If we choose homogeneous coordinates $x_i$ in $P(V^*)$ in
 which $Z$ consists of points $(1,0,0), (0,1,0), (0,0,1)$
then $c$ is given by the formula
$t_0 = x_1x_2, t_1 = x_0x_2, t_2 = x_0x_1$ where $t_i$ are appropriate
coordinates in $P(M^*)$.

c) The Schur quadric $Q\in P(M^*)$ is the conic $t_0^2 +t_1^2 + t_2^2 = 0$.

d) The curve $C(E)$ is the inverse image of this conic under
the Cremona
transformation defined in b). In other words, the equation of
 $C(E)$ is $x_0^2x_1^2 +
x_0^2x_2^2 + x_1^2x_2^2 = 0$.

e) The  subspaces $A_z$ are coordinate lines
$\{t_i = 0\}$ in $P(M^*)$,
the  subspaces $A'_z$
 are the opposite points of the coordinate
triangle i.e., points $\{t_i = t_j = 0\}$.

\vskip .3cm

\noindent {\bf 4.4. The case $n=4$.} The moduli space
 $M(-1,4)$ has dimension
12. As shown in [DK], an open dense subset in $M(-1,4)$
 is provided by
normalized logarithmic bundles
$$E_{\rm norm}({\cal H}) = \Omega^1_{P(V)}(\log {\cal H})
 \otimes {\cal O}(-2)$$
where ${\cal H} = (H_1,...,H_6)$ is an arrangement of 6
lines in $P(V) = P^2$
in general position. We consider only such bundles $E$.
 Let $p_i\in P(V^*)$
be points corresponding to lines $H_i \i P(V)$. We first
assume that $p_i$
do not lie on a conic (i.e., $H_i$ are not all tangent to
a conic). In this
case:

a) $Z = \{p_1,...,p_6\}$.

b) The variety $\Sigma \i P(M^*) = P^3$ is the cubic surface obtained
by blowing up $Z$.

c) The quadric $Q$ is the classical Schur quadric associated with the
double - six $\{l_i = A_{p_i},\,\, l'_i = A'_{p_i}\}$ (see \S 0).
This follows from Theorem 2.17.

d) The curve $C(E)$ is the image under $\pi_\Sigma: \Sigma \rightarrow
P(V^*)$ of the intersection $\Sigma\cap Q$. The intersection is non
singular of degree 6 and genus 4; the projection will have nodes at $p_i$
since  each the six lines $l_i \i \Sigma$ blown down to $p_i$ by
$\pi_\Sigma$
meets $Q$ twice.

e) The subspaces $A_{p_i}, A'_{p_i}$ form the standard double - six
associated
to the blow-down $\pi_\Sigma$.

\vskip .2cm

\noindent If all $p_i$ do lie on a conic $\Gamma \i P(V^*)$, the situation
changes. In this case $E({\cal H})$ is the Schwartzenberger bundle
associated to $\Gamma$ (see [Schw 1-2], [DK]) and
 the features are as follows:

a) $Z$ equals the conic $\Gamma$ (so dim $Z = 1$).

b) The variety $\Sigma$ is the union of a smooth quadric surface $Q$ and
a plane $\Pi$. The projection $\pi_\Sigma: \Sigma \rightarrow P(V^*)$
maps $\Pi$ bijectively to $P(V^*)$ and projects $Q = P^1 \times P^1$
to one of its $P^1$ - factors which is then being embedded into $P(V^*)$
as the conic $\Gamma$.

c) The "Schur quadric" is the surface $Q$ from n. b).

d) The curve $C(E)$ coincides with $\Gamma$ (taken three times).

e) For any $z\in Z = \Gamma$ the lines $A_z$ and
$A'_z$ both coincide with
the generator of $Q = P^1 \times P^1$ mapped into $z$
 by $\pi_\Sigma$, see n.b).

\vskip .3cm

\noindent {\bf 4.5. Bring's curve as $C(E)$.} Consider the situation of
Example 0.13: the cubic surface $\Sigma$ is given by equation
$x_1^3 + ... + x_5^3 = 0$ where $x_i$ are linear functions on $M^*$
constrained by $\sum x_i = 0$. The Schur quadric corresponding
to double - six described in n. 0.13 is given by $\sum x_i^2 = 0$.
The intersection $C = \Sigma\cap Q$ i.e. the curve given in $P^4$
by equations
$$ \sum x_i = \sum x_i^2 = \sum x_i^3 = 0$$
is known as  {\it Bring's curve} [K] [Hu 2].
The blow-down of the  first six lines of the double - six
described in n. 0.13
gives 6 points $p_1,...,p_6$
in  $P^2$ forming an orbit
 of the alternating group
$A_5$ [Hu 2]. These points will be the nodes of the sextic
curve $\pi_\Sigma (C) \i P^2$ i.e., of the
projection of $C$ to $P^2$, which is also called Bring's curve.
The equation of $\pi_\Sigma (C)$ can be found in [Hu 2], p. 82.

\vskip .1cm

Thus  Bring's curve can be represented as the curve of
 JLSK of a certain
bundle on $\check P^2$: the (normalized) logarithmic
 bundle of the configuration of lines dual to $p_i$.

\vskip .3cm

\noindent {\bf 4.6. Hulsbergen bundles.} Let $q_1,...,q_n$ be $n$ points
in general position in $P(V)$. There exists an $n$ -dimensional family
of stable rank 2 bundles $E = E(\lambda_1, ..., \lambda_n),
\lambda_i \in {\bf C}^*$
on $P(V)$ with $c_1 = -1, c_2 = n$ such that $\{q_1,...,q_n\}$ is the set
of zeros of a section of $E(1)$ (see [Hu 1]). They are called Hulsbergen
bundles. For such $E$ the subscheme $Z$ of jumping lines of $E$ is reduced
and consists of $n\choose 2$ lines $<q_i, q_j>$. We denote by $f_j$
linear functions on $V^*$ corresponding to $q_i \in P(V)$. The linear
system  of curves of degree $n-1$
through $Z$ has a basis
formed by the curves
$$F_j = \prod_{i\neq j} f_i = 0.$$
This  system maps the surface $S  = {\rm Bl}_Z(P(V^*))$ to the surface
$\Sigma \i P(M^*) = P^{n-1}$ given, in natural homogeneous
coordinates
$(t_1,...,t_{n})$, by equations
$$\biggl(\prod_{i=1}^n t_i\biggl) \biggl( \sum_{i=1}^n
{a_{ji}\over t_i}\biggl)
= 0, j=1,...,n-3$$
where $(a_{j1}, ..., a_{jn}), \, j=1,...,n-3$ is a
basis of the space of linear
relations among the vectors $f_i$. In the coordinates
 $t_i$ the "Schur quadric"
$Q$ is given by the equation $\sum c_i t_i^2 = 0$ so
the curve of JLSK
in $P(V^*)$ has the equation
$$\sum_{i=1}^n c_i F_i^2 = 0.$$
(cf. [Hu 1]. n. 10.5). Note that $p: S \rightarrow \Sigma$ blows down
the proper transforms of the lines $l_i$ to singular points
 of $\Sigma$ which have the coordinates $(1,0,...,0), ...,
(0,..., 1)$.  These are points
belong to $a_M^{-1}( {\rm Hom}\,\,(V^*, H^*)_1)$. So Hulsbergen bundles
are not $\Sigma$ - generic in the sense of n. 2.5, although they are
$Z$ - generic.

\beginsection

\centerline {\bf References.}

\vskip 1cm

\item{[ACGH]} E.Arbarello, M.Cornalba, Ph.Griffiths, J.Harris,
 Geometry of
algebraic curves, Vol.1 (Grund. Math. Wiss, {\bf 267}),
Springer-Verlag, 1990.

\item{[B]} H.Baker, Principles of geometry, Vol.3,
Cambridge Univ. Press,
1927.

\item{[Bu]} W.Burnside, On the double - six which admits a
group of 120 collineations into inself, {\it Proc. Camb.
Phil. Soc.}, {\bf 16} (1911), 418 - 420.

\item{[DG]} E.Davis, A.Geramita, Birational morphisms to $P^2$:
an ideal - theoretic perspective, {\it Math. Ann.}, {\bf 279}
(1988), 435 - 448.

\item{[DK]} I.Dolgachev, M.Kapranov, Arrangements of hyperplanes
and vector bundles on $P^n$, preprint 1992.

\item{[DO]} I.Dolgachev, D.Ortland, Finite point sets in
projective spaces and theta-functions, {\it Ast\'erisque},
 {\bf 167}, Soc. Math. France, 1989.

\item{[G]} A.Geramita, The curves seminar at Queen's, Vol.6, {\it Queens'
Papers in Pure and Appl. Math.}, {\bf 83} (1989).

\item{[Gi]} A.Gimigliano, On Veronesean surfaces, {\it Indag.
Math.}, {\bf 51} (1989), 71 - 85.

\item{H]} R.Hartshorne, Algebraic Geometry, Springer - Verlag, 1977.

\item{[Hu1]} K.Hulek, Stable rank 2 bundles on $P_2$ with $c_1$ odd,
{\it Math. Ann.}, {\bf 242} (1979), 241-266.

\item{[Hu2]} K.Hulek, Geometry of the Horroks-Mumford bundle, Proc. Symp.
Pure Math., Vol. 46, Part 2, p. 69 - 85, AMS, 1987.

\item{[Kap]} I.Kaplansky, Commutative rings, University of
Chicago Press, 1974.

\item{[K]} F.Klein, Lectures on the icosahedron and the solution of
equations of fifth degree, Dover Publ., New York, 1956.

\item {[L]} J. Le Potier, Fibr\'es stables de rang 2 sur $P_2({\bf C})$,
{\it Math. Ann.}, {\bf 241} (1979), 217 - 256.

\item {[M]} Y.Manin. Cubic forms, North-Holland, 1977.

\item{[OSS]} C.Okonek, M.Schneider, H.Spindler, Vector bundles on
complex projective spaces, (Progress in Math., Vol.3), Birkh\"auser,
 Boston 1980.

\item{[R]} T.G.Room, The geometry of determinantal loci, Cambridge Univ.
Press, 1937.

\item {[Schur]} F. Schur, \"Uber die durch collineare Grundgebilde
erzeugten Curven und Fl\"achen, {\it Math. Ann.}, {\bf 18} (1881), 1 -
32.

\item{[Schw1]} R.L.E.Schwarzenberger, Vector bundles on the projective
plane, {\it Proc. London Math. Soc.}, {\bf 11} (1961), 623 - 640.

\item{[Schw2]} R.L.E.Schwarzenberger, The secant bundle of a
 projective
variety, {\it Proc. London Math. Soc.}, {\bf 14} (1964),
369 - 384.

\vskip 2cm

Authors' addresses:

\noindent I.D.: Department of Mathematics, University of
Michigan, Ann Arbor MI 48109, email: IGOR.DOLGACHEV@um.cc.umich.edu

\vskip .3cm

\noindent M.K.: Department of Mathematics, Northwestern
 University, Evanston IL 60208, email: kapranov@chow.math.nwu.edu

\bye

 collineare Grundgebilde
erzeugten Curven und Fl\"achen, {\it Math. Ann.}, {\bf 18} (1881), 1 -
32.

\item{[Schw1]} R.L.E.Schwarzenberger, Vector bundles on the projective
plane, {\it Proc. London Math. Soc.}, {\bf 11} (1961), 623 - 640.

\item{[Schw2]} R.L.E.Schwarzenberger, The secant bundle of a
 projective
variety, {\it Proc. London Math. Soc.}, {\bf 14} (1964),
369 - 384.

\vskip 2cm

Authors' addresses:

\noindent I.D.: Department of Mathematics, University of
Michigan, Ann Arbor MI 48109, email: IGOR.DOLGACHEV@um.cc.umich.edu

\vskip .3cm

\noindent M.K.: Department of Mathematics, Northwestern
 University, Evanston IL 60208, email: kapranov@chow.math.nwu.edu

\bye